\documentclass[conference, 9pt]{IEEEtran}
\IEEEoverridecommandlockouts

\usepackage{cite}
\usepackage{amsmath,amssymb,amsfonts}
\usepackage{stmaryrd}
\usepackage{graphicx}
\usepackage{textcomp}
\usepackage{xcolor}

\usepackage{mathtools}
\usepackage{bm}
\usepackage{ascmac}
\usepackage{amsthm}
\usepackage{mathrsfs}
\usepackage{siunitx}
\usepackage{algpseudocode}
\usepackage{algorithm}
\usepackage{caption}
\usepackage{cite}

\newcommand{\R}{\mathbb{R}}
\newcommand{\N}{\mathbb{N}}

\newcommand{\argmin}[1][]{\underset{{#1}}{\mathrm{argmin}} \ }
\newcommand{\product}[3][]{\left\langle {#2} , {#3} \right\rangle_{#1}}

\newcommand{\vv}{\bm{v}}
\newcommand{\w}{\bm{w}}
\newcommand{\x}{\bm{x}}
\newcommand{\y}{\bm{y}}

\newcommand{\f}{\bm{f}}
\newcommand{\Li}{\bm{L}}

\newcommand{\G}{\bm{G}}

\newcommand{\A}{\bm{A}}
\newcommand{\B}{\bm{B}}

\newcommand{\Q}{\bm{Q}}
\newcommand{\TFBF}{\bm{T}_{\mathrm{FBF}}}
\newcommand{\Talpha}{\bm{T}_{\alpha}}
\newcommand{\uf}[1]{\bm{\mathfrak{f}}_{{#1}}^{\langle \mathrm{u} \rangle}}
\newcommand{\uG}{\bm{\mathfrak{G}}^{\langle \mathrm{u} \rangle}}
\newcommand{\st}{\mathrm{s.t.} \ }
\newcommand{\find}{\mathrm{find} \ }
\newcommand{\Id}{\mathrm{Id}}
\newcommand{\FBF}{\mathrm{FBF}}
\newcommand{\Fix}{\mathrm{Fix}}

\newcommand{\zer}{\mathrm{zer}}

\newcommand{\uV}{\bm{\mathcal{V}}^{\langle u \rangle}}

\theoremstyle{plain}
\newtheorem{dfn}{Definition}[section]
\newtheorem{asm}[dfn]{Assumption}
\newtheorem{prop}[dfn]{Proposition}
\newtheorem{lem}[dfn]{Lemma}
\newtheorem{thm}[dfn]{Theorem}

\theoremstyle{definition}
\newtheorem{prob}[dfn]{Problem}

\newtheorem{rem}[dfn]{Remark}

\renewcommand{\theenumi}{\roman{enumi}}
\renewcommand{\labelenumi}{(\theenumi)}

\renewcommand{\baselinestretch}{0.92}

\def\BibTeX{{\rm B\kern-.05em{\sc i\kern-.025em b}\kern-.08em
    T\kern-.1667em\lower.7ex\hbox{E}\kern-.125emX}}
\begin{document}

\title{
Hierarchical Nash Equilibrium over Variational Equilibria via \\
Fixed-point Set Expression of Quasi-nonexpansive Operator
\thanks{This work was supported by JSPS Grants-in-Aid (19H04134, 24K23885).}}

\author{\IEEEauthorblockN{Shota Matsuo, Keita Kume, and Isao Yamada}
\IEEEauthorblockA{\textit{Dept. of Information and Communications Engineering, Institute of Science Tokyo}\thanks{Tokyo Institute of Technology merged with Tokyo Medical and Dental University to form ``Institute of Science Tokyo'' on October 1, 2024.}, JAPAN\\
\{matsuo, kume, isao\}@sp.ict.e.titech.ac.jp}}

\maketitle

\begin{abstract}
The equilibrium selection problem in the generalized Nash equilibrium problem (GNEP) has recently been studied as an optimization problem, defined over the set of all variational equilibria achievable through a lower-level non-cooperative game among players. 
However, to make such a selection fair for every player, we have to rely on an unrealistic assumption, that is, the availability of a trusted center that does not induce any bias for every player. 
In this paper, we study a new equilibrium selection problem, named the hierarchical Nash equilibrium problem (HNEP), and propose an iterative algorithm for solving the HNEP. 
The HNEP is designed to ensure a fair selection without assuming any trusted center. 
More precisely, the HNEP is the GNEP for an upper-level non-cooperative game defined over the set of all variational equilibria of the lower-level non-cooperative game. 
The proposed algorithm for the HNEP is established by applying the hybrid steepest descent method to a variational inequality defined over the fixed point set of a quasi-nonexpansive operator. 
Numerical experiments show the effectiveness of the proposed equilibrium selection problem and its algorithmic solution. 
\end{abstract}


\setlength{\abovedisplayskip}{5.2pt} 
\setlength{\belowdisplayskip}{5.2pt} 

\section{Introduction}
\label{sec:intro}
Game theory dates back to the pioneering work of von Neumann and Morgenstern \cite{NM1944}, and the innovative idea of Nash equilibrium (NE), introduced by John Nash \cite{nash1950equilibrium,nash1951non}, triggered the drastic expansion of applications of game theory. 
The NE is a well-balanced solution of non-cooperative games, in which multiple decision-makers $\mathcal{I}:=\{1, \dots, m\}$, called players, aim to decrease respectively their cost functions as much as possible. 
At the NE, any player, say $i\in\mathcal{I}$, cannot decrease solely $i$'s cost function by changing $i$'s variable, called $i$'s strategy, as long as the other players' strategies are unchanged. 
The NE has been generalized \cite{facchineiGeneralizedNashEquilibrium2010,BC2021} and advanced \cite{luo1996mathematical, g.scutariEquilibriumSelectionPower2012, g.scutariRealComplexMonotone2014, e.benenatiOptimalSelectionTracking2023,w.heDistributedOptimalVariational2024a} toward one of ideal goals in a variety of modern engineering/social systems (see, e.g., \cite{10018595}), such as wireless communication systems \cite{a.deligiannisGameTheoreticPowerAllocation2017, w.wangNonConvexGeneralizedNash2022, y.chenQoEAwareDecentralizedTask2024}, smart grids systems \cite{i.atzeniNoncooperativeDayAheadBidding2014, zhengPeertopeerEnergyTrading2022}, and machine learning \cite{NIPS2014_5ca3e9b1, h.tembineDeepLearningMeets2020}. 
One of the unified definitions of the Nash equilibrium problem (NEP) can be stated in the following form:
\begin{prob}[{Generalized Nash equilibrium problem (GNEP) \cite{facchineiGeneralizedNashEquilibrium2010}, \cite[Exm. 3.4]{BC2021}}]
  \label{prob:GNEP}
  Consider a non-cooperative game among players in $\mathcal{I}:=\{1, \dots, m\}$. 
  We follow the notations used in \cite{BC2021} as:
  \begin{enumerate}
    \renewcommand{\labelenumi}{(\alph{enumi})}
    \item 
    For every $i \in \mathcal{I}$, the strategy of player $i$ is defined by $x_i \in \mathcal{H}_i$, where $\mathcal{H}_i$ is a finite-dimensional real Hilbert space. 
    \item
    Whole players' strategies and players' strategies other than player $i\in \mathcal{I}$ are denoted by $\x:=(x_1, \dots, x_m) \in \bm{\mathcal{H}}:=\bigtimes_{i \in \mathcal{I}}\mathcal{H}_i$ and $\x_{\smallsetminus i}:=(x_1, \dots, x_{i-1}, x_{i+1}, \dots, x_m)$ respectively. 
    \item
    For each $i \in \mathcal{I}$ and any $(x_i, \y) \in \mathcal{H}_i \times \bm{\mathcal{H}}$, we set $(x_i;\y_{\smallsetminus i})=(y_1, \dots, y_{i-1}, x_i, y_{i+1}, \dots, y_m)$. 
    \item 
    $\mathcal{G}$ is a finite-dimensional real Hilbert space. 
  \end{enumerate}
  Suppose that for each player $i \in \mathcal{I}$, his/her cost function is given by: 
  \begin{equation}
    \label{eq:cost}
    \x=(x_1, \dots, x_m) \mapsto \iota_{C_i}(x_i) + \f_i(\x) + (\iota_{D} \circ \Li)(\x), 
  \end{equation}
  where $C_i \subset \mathcal{H}_i$ and $D \subset \mathcal{G}$ are nonempty closed convex sets, 
  $\f_i:\bm{\mathcal{H}}\rightarrow\R$ satisfies that 
  $\f_i(\cdot;\x_{\smallsetminus i}):\mathcal{H}_i \rightarrow \mathbb{R}$ is convex and differentiable for every $\x \in \bm{\mathcal{H}}$, 
  and $\Li:\bm{\mathcal{H}}\rightarrow\mathcal{G}$ is a linear operator. 
  Then, the GNEP is given by
  \begin{align}
    \label{eq:GNEP_practice}
    \begin{split}
      &\find \x=(x_1, \dots, x_m)\in \bm{\mathcal{H}} \\
      &(\forall i \in \mathcal{I}) \ x_i \in \argmin \iota_{C_i}(\cdot) + \f_i(\cdot; \x_{\smallsetminus i}) + 
      \left( \iota_{D} \circ \Li \right)(\cdot ; \x_{\smallsetminus i}). 
    \end{split}
  \end{align}
  We call a solution of \eqref{eq:GNEP_practice} a generalized Nash equilibrium (GNE) and denote the solution set of \eqref{eq:GNEP_practice} by $\bm{\mathfrak{S}}_{\mathrm{GNE}}$.  
\end{prob}
\begin{figure}[t]
  \centering
  \includegraphics[scale = 0.35]{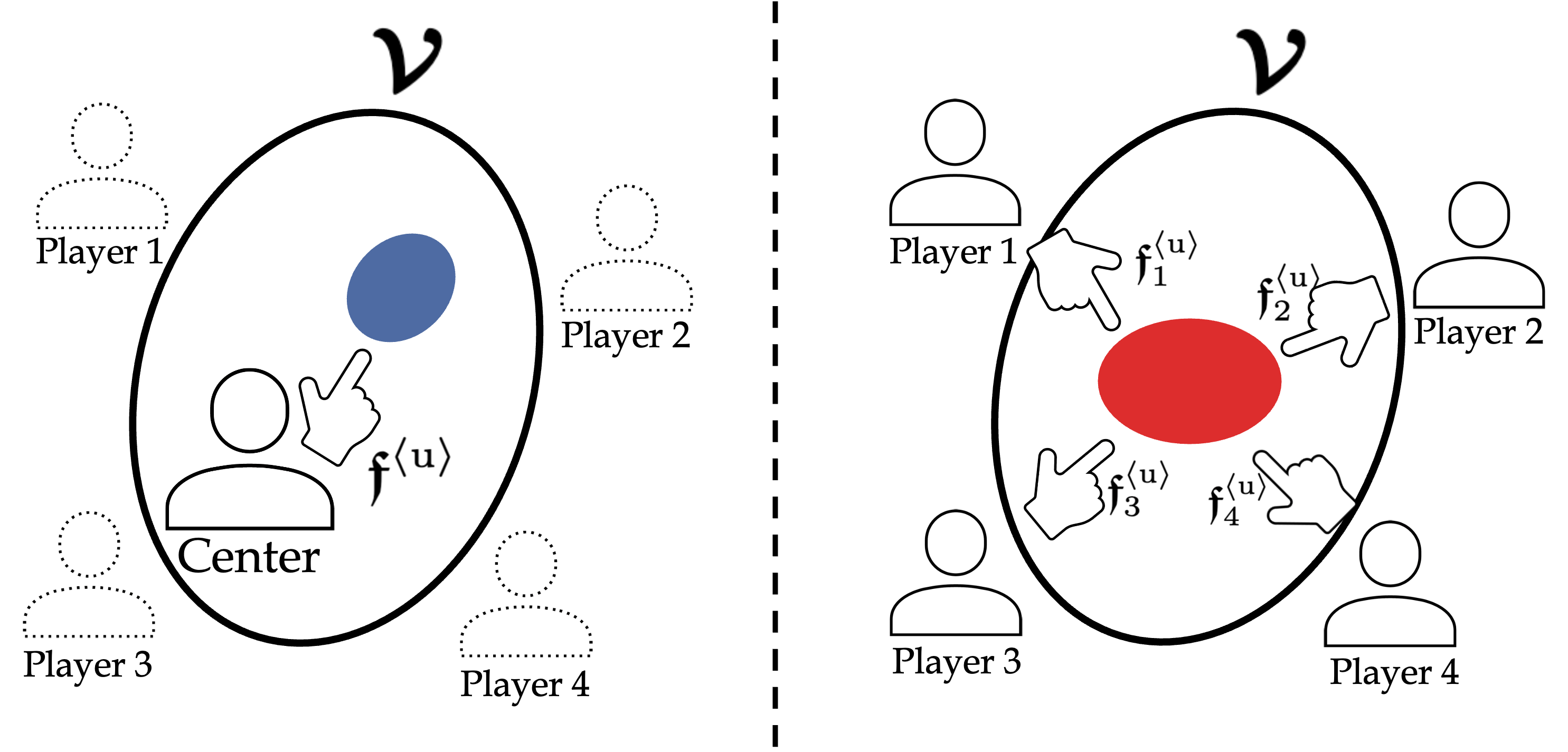}
  \caption{Conceptional comparison of two models for equilibrium selections (existing models $\llbracket$left$\rrbracket$ and proposed model $\llbracket$right$\rrbracket$) over $\bm{\mathcal{V}}$ of the lower-level non-cooperative game among all players $i \in \mathcal{I}$. \\
  \hspace*{1mm} $\llbracket$Left$\rrbracket$ The existing models \cite{g.scutariEquilibriumSelectionPower2012, g.scutariRealComplexMonotone2014, e.benenatiOptimalSelectionTracking2023,w.heDistributedOptimalVariational2024a} have been formulated to choose a special variational equilibrium by minimizing a single upper-level cost function $\uf{}$, designed hopefully by a trusted center. \\
  \hspace*{1mm} $\llbracket$Right$\rrbracket$ The proposed model is formulated to choose a special variational equilibrium, but in a different sense from the existing models, i.e., as an upper-level generalized Nash equilibrium (GNE) of a new non-cooperative game among all players $i\in\mathcal{I}$ with upper-level cost functions $\uf{i}$ designed by each player $i \in \mathcal{I}$.}
  \label{fig:image_HNEP}
\end{figure}

A specially valuable subset of $\bm{\mathfrak{S}}_{\mathrm{GNE}}$ is the solution set\footnote{For sufficient conditions to ensure $\bm{\mathcal{V}}\neq \varnothing$, see, e.g., \cite[Sec. 2.2]{FPvariational2003}, \cite{BC2021}.}
\begin{align}
  \label{eq:VE}
  \bm{\mathcal{V}}:=\left\{\vv \in \bm{\mathfrak{C}} \mid (\forall \w \in \bm{\mathfrak{C}}) \ \product[\bm{\mathcal{H}}]{\G(\vv)}{\w-\vv} \geq 0 \right\} \subset\bm{\mathcal{H}}
\end{align}
of the variational inequality $\mathrm{VI}(\bm{\mathfrak{C}}, \G)$ \cite[Definition 1.1.1]{FPvariational2003}, \cite{BC2017}, 
where 
$\bm{\mathfrak{C}}:=\{\x \in \bigtimes_{i \in \mathcal{I}} C_i \mid \Li\x \in D\}$ and 
\begin{equation}
  \label{eq:G}
  \bm{G}:\bm{\mathcal{H}}\rightarrow\bm{\mathcal{H}}:\x\mapsto(\nabla_1 \f_1(\x), \cdots, \nabla_m \f_m(\x))
\end{equation}
is defined with gradients $\nabla_i \f_i(\cdot;\x_{\smallsetminus i}):\mathcal{H}_i\rightarrow\mathcal{H}_i$ of $\f_i(\cdot;\x_{\smallsetminus i})$. 
Indeed, a point in $\bm{\mathcal{V}}$ is referred to as a {\em variational equilibrium} \cite{rosenExistenceUniquenessEquilibrium1965, facchineiGeneralizedNashGames2007}, \cite[Def. 3]{facchineiGeneralizedNashEquilibrium2010} and known to enjoy several desirable properties, such as fairness and larger social stability than any GNE in $\bm{\mathfrak{S}}_{\mathrm{GNE}}\setminus \bm{\mathcal{V}}$ \cite{facchineiGeneralizedNashEquilibrium2010, kulkarniVariationalEquilibriumRefinement2012}. 
Recent applications of variational equilibrium are found, e.g., in distributed control and signal processing over networks 
\cite{yiOperatorSplittingApproach2019, g.belgioiosoDistributedGeneralizedNash2022, l.ranDistributedGeneralizedNash2024}. 

In general, the set $\bm{\mathcal{V}}$ in \eqref{eq:VE} is an infinite set. 
This situation induces a naive question: can we design a fair mechanism for each player to reach a certainly desirable variational equilibrium in $\bm{\mathcal{V}}$ ? 
Regarding this question, \cite{g.scutariEquilibriumSelectionPower2012, g.scutariRealComplexMonotone2014, e.benenatiOptimalSelectionTracking2023, w.heDistributedOptimalVariational2024a} proposed to formulate a hierarchical convex optimization problem (see, e.g., \cite{Yamada-Yamagishi19}), i.e., minimization of a single {\em upper-level} convex function, say $\uf{}:\bm{\mathcal{H}}\to\R$, over $\bm{\mathcal{V}}$ in \eqref{eq:VE}, i.e., the set of all variational equilibria of the ({\em lower-level}) non-cooperative game. 
Indeed, \cite{g.scutariEquilibriumSelectionPower2012, g.scutariRealComplexMonotone2014} 
proposed iterative algorithms of nested structures by introducing an inner loop to solve certain subproblems. 
Quite recently, \cite{e.benenatiOptimalSelectionTracking2023, w.heDistributedOptimalVariational2024a} 
proposed to apply the hybrid steepest descent method 
\cite{Y2001, oguraNonstrictlyConvexMinimization2003, yamadaHybridSteepestDescent2005,Yamada-Yamagishi19} (see also \cite[Prop. 42]{p.l.combettesFixedPointStrategies2021}) 
to such hierarchical convex optimization problems without any inner loop. 

From the viewpoint of non-cooperative game theory, which has been pursuing ideal fairness among multiple players $\mathcal{I}$ without requiring any intervenient, say \emph{center} in this paper, another question arises: who in the world can design such a function $\uf{}$ certainly according to each player's hope without causing any risk of unexpected bias among players ?
If we assumed the availability of a center perfectly reliable to all players, we could delegate the authority of designing of $\uf{}$ to the center (we do not believe the availability of such a trusted center).

In this paper, by revisiting the spirit of John Nash, we resolve this dilemma without requiring either such a trusted center (see Fig. \ref{fig:image_HNEP}) or any randomness assumption\footnote{For example, random selection from possibly infinitely many VEs.}. 
More precisely, we formulate the {\em{hierarchical Nash equilibrium problem} (HNEP)} below, as a novel equilibrium selection problem, and propose an iterative algorithm to find a solution of the HNEP\footnote{$\iota_{C_i}$ and $\iota_{D}$ can be relaxed to more general prox-friendly functions in $\Gamma_0(\mathcal{H}_i)$ and $\Gamma_0(\mathcal{G})$, respectively.}.
\begin{prob}[Hierarchical Nash equilibrium problem (HNEP)]
  \label{prob:HierarchicalGame_0}
  Under the setting of the lower-level non-cooperative game formulated in the form of Problem \ref{prob:GNEP}, let $\bm{\mathcal{V}}\neq \varnothing$ in \eqref{eq:VE} be the set of all variational equilibria of the lower-level game. Then the HNEP is given, as an GNEP of an \emph{upper-level} non-cooperative game over $\bm{\mathcal{V}}$, by
  \begin{align}
    \label{eq:HierarchicalGame_0}
    \begin{split}
      \find &\x = (x_1, \dots, x_m) \in \bm{\mathcal{H}} \\
      \st &(\forall i \in \mathcal{I}) \ x_i \in \argmin \uf{i}(\cdot; \x_{\smallsetminus i}) + 
      \iota_{\bm{\mathcal{V}}}(\cdot; \x_{\smallsetminus i}), 
    \end{split}
  \end{align}
  where $\uf{i}:\bm{\mathcal{H}}\rightarrow\R$ is the player $i$'s {\em upper-level} cost function such that $\uf{i}(\cdot;\x_{\smallsetminus i}):\mathcal{H}_i \rightarrow \mathbb{R}$ is convex and differentiable for every $\x \in \bm{\mathcal{H}}$ (see Remark \ref{rem:f_u} (i) regarding the design of $\uf{i}$).
\end{prob}
\begin{rem}[On the HNEP \eqref{eq:HierarchicalGame_0}]
  \label{rem:f_u} \ 
  \begin{enumerate}
    \item $\uf{i}:\bm{\mathcal{H}}\rightarrow\R$ can be designed, by each player $i(\in\mathcal{I})$, independently of the lower-level non-cooperative game. 
    \item By assigning a common function $\uf{}$ to $\uf{i} \ (\forall i \in \mathcal{I})$, the HNEP \eqref{eq:HierarchicalGame_0} is reduced to the hierarchical convex optimization problems over $\bm{\mathcal{V}}$ \cite{g.scutariEquilibriumSelectionPower2012, g.scutariRealComplexMonotone2014, e.benenatiOptimalSelectionTracking2023,w.heDistributedOptimalVariational2024a}. 
  \end{enumerate}
\end{rem}
For finding a solution of Problem \ref{prob:HierarchicalGame_0}, the proposed algorithm (Algorithm 1 in Sec. \ref{sec:HVE}) is established by applying the hybrid steepest descent method \cite{yamadaHybridSteepestDescent2005} to a fixed point expression of $\bm{\mathcal{V}}$ with a quasi-nonexpansive operator \cite{tsengModifiedForwardBackwardSplitting2000}. 
In numerical experiments (Section \ref{sec:expriment}), we apply Algorithm 1 to the HNEP in a scenario of the linearly-coupled aggregative game \cite[Section IV.B]{g.belgioiosoSemiDecentralizedGeneralizedNash2023a}. The numerical result demonstrates that Algorithm 1 achieves a variational equilibrium of the upper-level non-cooperative game. 

\textbf{Notation}:
We denote the sets of all nonnegative integers, real numbers, and positive real numbers by $\N$, $\R$, and $\R_{++}$, respectively. 
Let $(\mathcal{H}, \product[\mathcal{H}]{\cdot}{\cdot}, \|\cdot\|_{\mathcal{H}})$ and 
$(\mathcal{K}, \product[\mathcal{K}]{\cdot}{\cdot}, \|\cdot\|_{\mathcal{K}})$ 
be finite-dimensional real Hilbert spaces. 
The notation $\bm{0}_{\mathcal{H}}$ is the zero element of $\mathcal{H}$. 
For a linear operator $L:\mathcal{H}\to\mathcal{K}$, 
the operator norm of $L$ is $\|L\|_{\mathrm{op}}:=\sup_{\|x\|_{\mathcal{H}} \leq 1}\|Lx\|_{\mathcal{K}}$ 
and the adjoint of $L$ is the unique operator $L^*:\mathcal{K}\to\mathcal{H}$ that satisfies 
$(\forall x\in\mathcal{H}, \forall y \in \mathcal{K}) \ \langle Lx, y\rangle_{\mathcal{K}} = \langle x, L^*y\rangle_{\mathcal{H}}$. 
The identity operator is denoted by $\Id$. 
The fixed point set of an operator $T:\mathcal{H}\rightarrow\mathcal{H}$ is 
$\Fix(T) := \{x\in\mathcal{H}\mid T(x)=x\}$. 
The set of zeros of a set-valued operator $A:\mathcal{H}\rightarrow 2^{\mathcal{H}}$ 
is $\zer(A) := \{x\in\mathcal{H}\mid \bm{0}_{\mathcal{H}}\in A(x)\}$. 
A function $f:\mathcal{H}\rightarrow (0, \infty]$ is called 
(i) proper if $\mathrm{dom}(f) := \{x\in\mathcal{H}\mid f(x)<\infty\} \neq \varnothing$, 
(ii) lower semicontinuous if $\{x\in\mathcal{H}\mid f(x) \leq \alpha\}$ is closed 
for every $\alpha \in \R$, 
(iii) convex if $(\forall x, y \in \mathcal{H}) \ (\forall \lambda \in [0, 1]) \ 
f(\lambda x + (1-\lambda)y) \leq \lambda f(x) + (1-\lambda)f(y)$. 
We denote the set of proper, lower semicontinuous, and convex 
functions on $\mathcal{H}$ by $\Gamma_0(\mathcal{H})$. 
For $f\in \Gamma_0(\mathcal{H})$, 
the conjugate of $f$ is $f^*:\mathcal{H}\to [-\infty, \infty]:u\mapsto\sup_{x\in\mathcal{H}}(\langle x, u\rangle_{\mathcal{H}} - f(x))$, 
the subdifferential of $f$ at $x\in\mathcal{H}$ is 
$\partial f(x) := \{v\in\mathcal{H}\mid 
(\forall y \in \mathcal{H}) \ f(y) \geq f(x) + \product[\mathcal{H}]{\cdot}{v}\}$. 
Let $C \subset \mathcal{H}$ be a nonempty closed convex set. 
The indicator function of $C$ is 
$\iota_{C}(x) = 0$ if $x \in C$ and $\iota_{C}(x) = \infty$ otherwise. 
The normal cone operator to $C$ at $x\in\mathcal{H}$ is 
$N_C(x) = \{v\in\mathcal{H}\mid 
(\forall y \in C) \ \product[\mathcal{H}]{v}{y-x} \leq 0\}$ if $x\in C$ and $N_{C}(x) = \varnothing$ otherwise.
The projection mapping onto $C$ is 
$P_C:\mathcal{H}\rightarrow\mathcal{H}:x\mapsto \textstyle\mathrm{argmin}_{y\in C}\|y - x\|_{\mathcal{H}}$. 
\section{Preliminaries}
\label{sec:preliminaries}
Fixed point theory provides powerful tools for solving the GNEP \cite{p.l.combettesFixedPointStrategies2021, g.belgioiosoDistributedGeneralizedNash2022}. 
In this section, we present a fixed point expression of $\bm{\mathcal{V}}$ in \eqref{eq:VE} under the following standard assumption \cite{yiOperatorSplittingApproach2019,e.benenatiOptimalSelectionTracking2023,BC2021}. 
\begin{asm}
  \label{asm:monotone}
  Under the setting of Problem \ref{prob:GNEP}, assume that 
  $\G:\bm{\mathcal{H}}\rightarrow\bm{\mathcal{H}}$ in \eqref{eq:G}
  is monotone and $\kappa_{\G}$-Lipschitzian, i.e., 
  \begin{align}
    &(\forall \x, \y \in \bm{\mathcal{H}}) \ 
    \product[\bm{\mathcal{H}}]{\G(\x) - \G(\y)}{\x - \y} \geq 0, \\
    &(\exists \kappa_{\G} > 0) \ (\forall \x, \y \in \bm{\mathcal{H}}) \ 
    \|\G(\x) - \G(\y)\|_{\bm{\mathcal{H}}} \leq \kappa_{\G} \|\x - \y\|_{\bm{\mathcal{H}}}. 
  \end{align}
\end{asm}
\begin{prop}[{Fixed point expression of $\bm{\mathcal{V}}$}]
  \label{prop:zero_closed_convex}
  Under $\bm{\mathcal{V}}\neq\varnothing$ and Assumption \ref{asm:monotone}, $\bm{\mathcal{V}}$ is closed convex. 
  In addition, by defining  
  $\TFBF:\bm{\mathcal{H}} \times\mathcal{G} \rightarrow \bm{\mathcal{H}} \times\mathcal{G}$, 
  with $\gamma \in (0, 1/(\kappa_{\G} + \|\Li\|_{\mathrm{op}}))$, as 
  \begin{equation}
    \label{eq:FBF}
    \TFBF := (\Id - \gamma\A) \circ (\Id + \gamma \B)^{-1} \circ (\Id - \gamma \A) + \gamma \A, 
  \end{equation}
  and its $\alpha(\in(0, 1))$-averaged operator 
  \begin{equation}
    \label{eq:Talpha}
    \Talpha := (1-\alpha)\Id + \alpha \TFBF, 
  \end{equation}
  where
  \begin{align}
    \label{eq:A}
    &\A:\bm{\mathcal{H}} \times\mathcal{G} \rightarrow \bm{\mathcal{H}} \times\mathcal{G} : (\x, u) \mapsto (\G(\x)+\Li^*u, -\Li\x), \\
    \label{eq:B}
    &\B:\bm{\mathcal{H}} \times\mathcal{G} \rightarrow 2^{\bm{\mathcal{H}}}\times 2^{\mathcal{G}} : (\x, u) \mapsto \big(\textstyle\bigtimes_{i \in \mathcal{I}} N_{C_i}(x_i) \big) \times 
    \partial \iota_{D}^*(u), 
  \end{align}
  we have $\Fix(\TFBF) = \Fix(\Talpha)$ and 
  \begin{enumerate}
    \item 
    $\bm{\mathcal{V}}$ can be expressed as $\bm{\mathcal{V}} = \Q_{\bm{\mathcal{H}}}\bigl( \Fix(\Talpha) \bigr) \neq \varnothing$ 
    with a canonical projection $\Q_{\bm{\mathcal{H}}}:\bm{\mathcal{H}} \times\mathcal{G} 
    \rightarrow \bm{\mathcal{H}}:(\x, u)\mapsto\x$~onto~$\bm{\mathcal{H}}$~\cite{e.benenatiOptimalSelectionTracking2023}. 
    \item $\Talpha$ is strongly attracting quasi-nonexpansive\!\!
    \footnote{
      The condition \eqref{eq:strongly_attracting} automatically implies that $\Talpha$ is \emph{quasi-nonexpansive}, i.e., 
      $(\forall \bm{\xi}:=(\x, u) \in \bm{\mathcal{H}}\times\mathcal{G}, \forall \bm{\zeta} \in \Fix(\Talpha)) \ \|\Talpha(\bm{\xi}) - \bm{\zeta}\| \leq \|\bm{\xi} - \bm{\zeta}\|$, 
      and thus $\Fix(\Talpha) \subset \bm{\mathcal{H}}\times\mathcal{G}$ is closed convex. 
    }\cite{yamadaHybridSteepestDescent2005}, i.e., 
    \begin{align}
      \label{eq:strongly_attracting}
      \begin{split}
        &(\forall \bm{\xi}:=(\x, u) \in \bm{\mathcal{H}}\times\mathcal{G}, \forall \bm{\zeta} \in \Fix(\Talpha)) \\
        &\qquad \frac{1-\alpha}{\alpha}\|\Talpha(\bm{\xi}) - \bm{\xi}\|^2  \leq \|\bm{\xi} - \bm{\zeta}\|^2 - \|\Talpha(\bm{\xi}) - \bm{\zeta}\|^2. 
      \end{split}
    \end{align}
    Moreover, $\Talpha$ is \emph{quasi-shrinking}\!\!
    \footnote{
      Due to the space limitation, we omit the detailed explanation on quasi-shrinkingness. 
      See \cite[Sec. 3]{yamadaHybridSteepestDescent2005}, \cite[Prop. 2.10]{cegielskiAlgorithmSolvingVariational2013} for details. 
    }on any bounded, closed, and convex set 
    $\bm{\mathcal{C}}\subset\bm{\mathcal{H}}\times\mathcal{G}$ satisfying $\Fix(\Talpha)\cap\bm{\mathcal{C}} \neq \varnothing$. 
  \end{enumerate}
\end{prop}
\section{Proposed algorithm for solving HNEP}
\label{sec:HVE}
Recall that the set $\bm{\mathcal{V}}$, in \eqref{eq:VE}, of all variational equilibria is a valuable subset of the solution set of the GNEP \eqref{eq:GNEP_practice} in the sense of fairness. 
In the following, to apply this remarkable feature achievable by variational equilibrium to  equilibrium selection from $\bm{\mathcal{V}}$, we focus on finding a point in  
\begin{equation}
  \label{eq:VI_HierarchicalGame}
  \uV := \biggl\{\x \in \bm{\mathcal{V}} \mid 
  \product[\bm{\mathcal{H}}]{\uG(\x)}{\y -\x} \geq 0 \ (\forall \y \in \bm{\mathcal{V}})\biggr\}, 
\end{equation}
where 
\begin{equation}
  \label{eq:uG}
  \uG : \bm{\mathcal{H}}\rightarrow\bm{\mathcal{H}} : \x \mapsto \left( \nabla_1 \uf{1}(\x), \dots, \nabla_m \uf{m}(\x) \right) 
\end{equation}
with gradients $\nabla_i \uf{i}(\cdot;\x_{\smallsetminus i})$ of $\uf{i}(\cdot;\x_{\smallsetminus i})$ in Problem \ref{prob:HierarchicalGame_0}. 
Indeed, by letting $(C_i, D, \G, \Li):=(\mathcal{H}_i, \bm{\mathcal{V}}, \uG, \Id)$ in \eqref{eq:VE} and the closed convexity of $\bm{\mathcal{V}}$ (see Proposition \ref{prop:zero_closed_convex}), we see that
\begin{enumerate}
  \item $\uV$ is the solution set of $\mathrm{VI}(\bm{\mathcal{V}}, \uG)$, 
  \item $\uV$ is a valuable subset of the solution set of Problem \ref{prob:HierarchicalGame_0}. 
\end{enumerate}
Note that, to $\mathrm{VI}(\bm{\mathcal{V}}, \uG)$, we can not directly apply standard projected gradient type algorithms (e.g., \cite{goldstein1964convex}) because the projection onto $\bm{\mathcal{V}}$ is not available. 
To overcome this difficulty, we use a translation of $\mathrm{VI}(\bm{\mathcal{V}}, \uG)$ into a variational inequality over the fixed point set of $\Talpha$ in \eqref{eq:Talpha}. 
\begin{lem}
  Under the settings of Problem \ref{prob:HierarchicalGame_0} and Proposition \ref{prop:zero_closed_convex}, $\uV$ can be expressed as 
  \begin{equation}
    \label{eq:VE_u_F}
    \uV = \Q_{\bm{\mathcal{H}}}(\uV_{\FBF}), 
  \end{equation}
  where $\Q_{\bm{\mathcal{H}}}:\bm{\mathcal{H}}\times\mathcal{G} \rightarrow \bm{\mathcal{H}}$ is the 
  canonical projection onto $\bm{\mathcal{H}}$,
  $\uV_{\FBF} \subset \bm{\mathcal{H}}\times\mathcal{G}$ is the set of all solutions~of~$\mathrm{VI}(\Fix(\Talpha), \uG_{\FBF})$,~i.e., 
  \begin{align}
    \begin{split}
      \label{eq:VI_HierarchicalGame_F}
      &\find \bm{\xi}^{\star} \in \Fix(\Talpha) \ \st \\
      &(\forall \bm{\zeta} \in \Fix(\Talpha)) \ \product[]{\uG_{\FBF}(\bm{\xi}^{\star})}{\bm{\zeta} - \bm{\xi}^{\star}}\geq 0, 
    \end{split}
  \end{align}
  and $\uG_{\FBF}:\bm{\mathcal{H}}\times\mathcal{G}\rightarrow\bm{\mathcal{H}}\times\mathcal{G}:\bm{\xi}:=(\x, u) \mapsto (\uG(\x), \bm{0}_{\mathcal{G}})$. 
\end{lem}
By noting that the quasi-nonexpansive mapping $\Talpha$ enjoys the quasi-shrinking condition, we can apply 
the hybrid steepest descent method\!\!
\footnote{
  The hybrid steepest descent method \cite{yamadaHybridSteepestDescent2005} 
  has been proposed for the variational inequality problem over the 
  fixed point set of a quasi-nonexpansive mapping. 
  Under the quasi-shrinking condition on the mapping, 
  the method has a guarantee of convergence to a desired solution of such a variational inequality problem 
  \cite[Theorem 5]{yamadaHybridSteepestDescent2005}. 
}\cite[Theorem 5]{yamadaHybridSteepestDescent2005} to the problem \eqref{eq:VI_HierarchicalGame_F}:
\begin{equation}
  \label{eq:HSDM}
  (n\in\N) \ \bm{\xi}_{n+1} = \Talpha(\bm{\xi}_n) - \lambda_{n+1} \uG_{\FBF}(\Talpha(\bm{\xi}_n)) 
\end{equation}
with an arbitrarily given initial point $\bm{\xi}_0 \in \bm{\mathcal{H}}\times\mathcal{G}$ and 
stepsize $(\lambda)_{n\in\N} \subset [0, \infty)$ satisfying (H1) $\lim_{n\rightarrow \infty} \lambda_n = 0$ 
and (H2) $\sum_{n \in \N} \lambda_n = \infty$, e.g., $\lambda_n = 1/n$. 
Algorithm \ref{alg:HSDM} illustrates a concrete expression of \eqref{eq:HSDM}. 
Regarding the convergence of Algorithm \ref{alg:HSDM}, the following theorem is available under Assumption \ref{asm:paramonotone} below. 

\begin{algorithm}[tbp]
  \caption{Hybrid steepest descent method for \eqref{eq:VI_HierarchicalGame_F}}
  \label{alg:HSDM}
  \begin{algorithmic}[1]
    \State $\textbf{Input}: \gamma \in \left(0, 1/(\kappa_{\G} + \|\Li\|_{\mathrm{op}}) \right), 
    \alpha \in (0, 1), (\lambda_n)_{n \in \N} \subset [0, \infty) \ \textrm{satisfying} \ 
    \textrm{(H1)} \ \textrm{and} \ \textrm{(H2)}, (\x_0, u_0) = 
    ((x_{i, 0})_{i \in \mathcal{I}}, u_{0}) \in \bm{\mathcal{H}}\times\mathcal{G}. $
    \State $\textbf{Set}: \Pi_i:\bm{\mathcal{H}}\rightarrow\mathcal{H}_i:\x\mapsto x_i. $
    \State $\textbf{for } n = 1, 2, \dots \textbf{ do}$
      \State $\quad \textrm{Forward-backward} \ \textrm{step}:$
      \State $\qquad (\forall i \in \mathcal{I}) \ y_{i, n} \gets P_{C_i} \bigl[x_{i, n} - \gamma \bigl(\nabla_i \f_i(\x_n) + \Pi_i (\Li^*u_n) \bigr) \bigr]$
      \State $\qquad w_n \gets u_n - \gamma P_{D}\bigl[(1/\gamma)u_{n} + \Li \x_n \bigr]$
      \State $\quad \textrm{Forward} \ \textrm{step}:$
      \State $\qquad (\forall i \in \mathcal{I}) \ \tilde{y}_{i, n} \gets y_{i, n} - \gamma \bigl[ \bigl(\nabla_i \f_i(\y_n) + \Pi_i (\Li^* w_n)\bigr) - $
      \State $\qquad \qquad \qquad \qquad \ \ \bigl(\nabla_i \f_i(\x_n) + \Pi_i (\Li^* u_n)\bigr) \bigr]$
      \State $\qquad \tilde{w}_{n} \gets w_{n} + \gamma (\Li \y_n - \Li \x_n)$
      \State $\quad \alpha\textrm{-averaging} \ \textrm{step}:$
      \State $\qquad (\x_{n+1/2}, u_{n+1}) \gets (1-\alpha)(\x_n, u_n) + \alpha(\tilde{\y}_n, \tilde{w}_n)$
      \State $\quad \textrm{Steepest} \ \textrm{descent} \ \textrm{step}:$
      \State $\qquad (\forall i \in \mathcal{I}) \ x_{i, n+1} \gets x_{i, n+1/2} - \lambda_{n+1} \nabla_i \uf{i}(\x_{n+1/2})$
    \State $\textbf{end for}$
  \end{algorithmic}
\end{algorithm}

\begin{asm}[Assumption for convergence of \eqref{eq:HSDM}]
  \label{asm:paramonotone}
  Under the settings of Problem \ref{prob:HierarchicalGame_0} and Proposition \ref{prop:zero_closed_convex}, assume that 
  \begin{enumerate}
    \item $\uG:\bm{\mathcal{H}}\rightarrow\bm{\mathcal{H}}$ defined by \eqref{eq:uG} is 
          Lipschitzian and paramonotone, 
          i.e., $\uG$ is monotone and 
          for any $\x, \y \in \bm{\mathcal{H}}$, the following holds: 
          \begin{equation}
            \bigl\langle \uG(\x) - \uG(\y),  \x - \y \bigr\rangle_{\bm{\mathcal{H}}} = 0
            \Leftrightarrow \uG(\x) = \uG(\y). 
          \end{equation}
    \item $\Fix(\Talpha)$ is bounded. 
    \item $(\bm{\xi}_n)_{n\in\N} = (\x_n, u_n)_{n\in\N}$ generated by \eqref{eq:HSDM} (Algorithm \ref{alg:HSDM}) is bounded. 
  \end{enumerate}
\end{asm}
\begin{thm}[Convergence of Algorithm \ref{alg:HSDM}]
  \label{thm:HSDM_convergence}
  Under Assumption \ref{asm:paramonotone},  
  we define $\uV$ as \eqref{eq:VI_HierarchicalGame} and the distance 
  $d(\cdot, \uV):\bm{\mathcal{H}} \ni \x \mapsto \min_{\y \in \uV}\| \y - \x \|_{\bm{\mathcal{H}}}$ 
  to the nonempty\!\!
  \footnote{
    Under Assumption \ref{asm:paramonotone}, $\uV \neq \varnothing$ 
    is automatically guaranteed (see \cite[Theorem 5]{yamadaHybridSteepestDescent2005}). 
  }closed convex set $\uV$. 
  Then, for any initial point $\bm{\xi}_0 \in \bm{\mathcal{H}}\times\mathcal{G}$, 
  the sequence $(\bm{\xi}_n)_{n \in \N} = (\x_n, u_n)_{n\in\N}$ generated by \eqref{eq:HSDM} (Algorithm \ref{alg:HSDM})
  enjoys $\lim_{n \rightarrow \infty} d(\x_n, \uV) = 0$. 
\end{thm}

\begin{rem}[Elimination of possible concerns regarding Assumption \ref{asm:paramonotone}]
  \label{rem:ass} \ 
  \begin{enumerate}
    \item Regarding Assumption \ref{asm:paramonotone} (i), 
          many non-cooperative games are known to satisfy this condition, e.g., 
          the Nash-Cournot game, 
          which is appeared in the control of electric vehicles \cite{yiOperatorSplittingApproach2019}. 
          Moreover, if $\uG$ is strictly monotone, i.e., for any $\x, \y \in \bm{\mathcal{H}}$, 
          \begin{equation}
            \x \neq \y \Rightarrow 
            \bigl\langle \uG(\x) - \uG(\y),  \x - \y \bigr\rangle_{\bm{\mathcal{H}}} > 0, 
          \end{equation}
          then $\uG$ is paramonotone \cite[Chap. 22]{BC2017}. 
    \item Regarding Assumption \ref{asm:paramonotone} (ii) and (iii), 
          the boundedness of $\Fix(\Talpha)$ and that of $(\bm{\xi})_{n \in \N}$ are not automatically guaranteed. 
          However, these conditions may not be restrictive for most practitioners by just modifying 
          our original target \eqref{eq:VI_HierarchicalGame_F} into 
          \begin{align}
            \begin{split}
              \label{eq:VI_HierarchicalGame_F_B}
              &\find \bm{\xi}^{\star} \in \overline{B}(0, r)\cap\Fix(\Talpha) \ \st \\
              &(\forall \bm{\zeta} \in \overline{B}(0, r)\cap\Fix(\Talpha)) \ \product[]{\uG_{\FBF}(\bm{\xi}^{\star})}{\bm{\zeta} - \bm{\xi}^{\star}}\geq 0, 
            \end{split}
          \end{align}
          with a sufficiently large closed ball $\overline{B}(0, r) \subset \bm{\mathcal{H}}\times\mathcal{G}$. 
          Note that the hybrid steepest descent method \cite[Theorem 5]{yamadaHybridSteepestDescent2005} is 
          applicable to \eqref{eq:VI_HierarchicalGame_F_B} because $P_{\overline{B}(0, r)} \circ \Talpha$ is a 
          strongly attracting quasi-nonexpansive mapping enjoying the quasi-shrinking condition and 
          $\Fix(P_{\overline{B}(0, r)} \circ \Talpha) =\overline{B}(0, r)\cap\Fix(\Talpha)$. 
  \end{enumerate}
\end{rem}
\section{Application of the hierarchical Nash equilibrium problem to linearly-coupled aggregative game}
\label{sec:expriment}
\begin{figure*}
  \centering
  \includegraphics[scale=0.48]{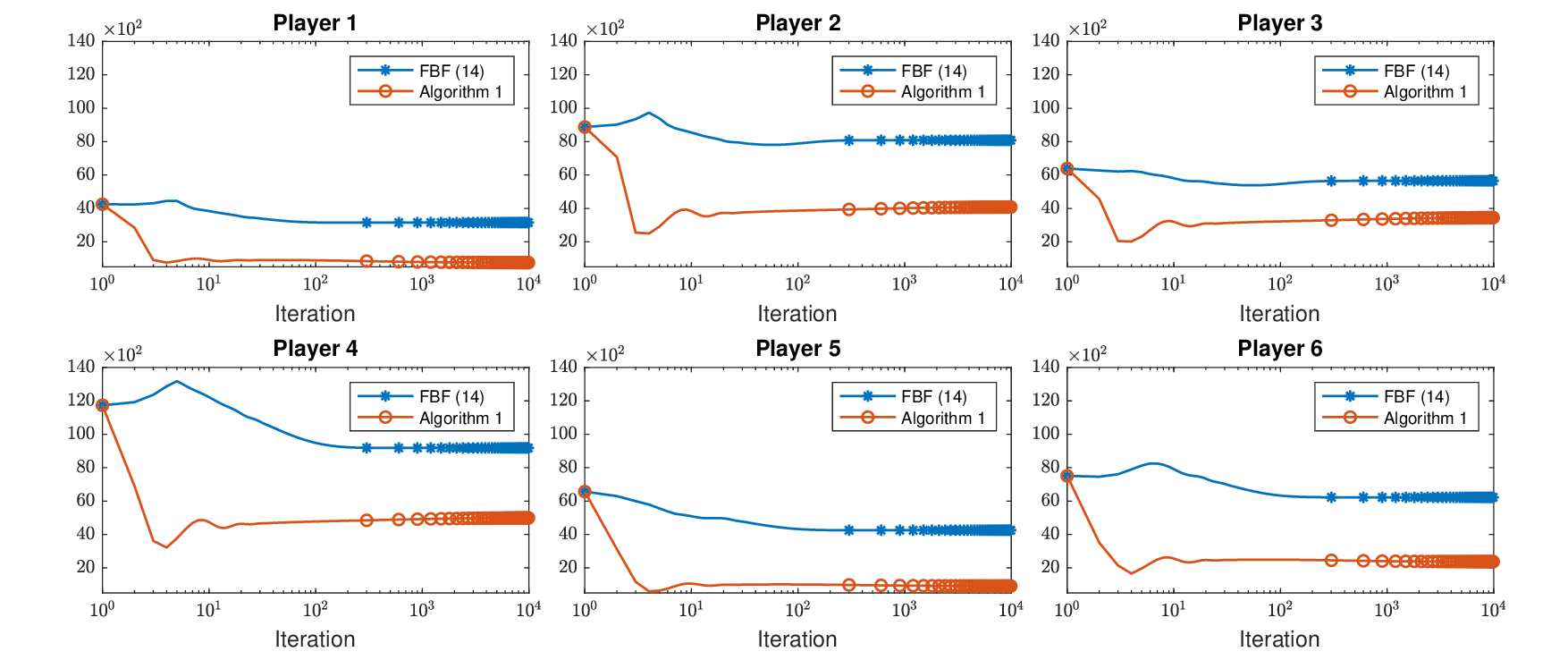}
  \vspace*{-0.15cm}
  \renewcommand{\thefigure}{3}
  \caption{Values of each player's upper-level cost function $\uf{i} \ (i\in\mathcal{I})$}
  \label{fig:cost_u}
\end{figure*}
\begin{figure}[t]
  \centering
  \includegraphics[scale = 0.53]{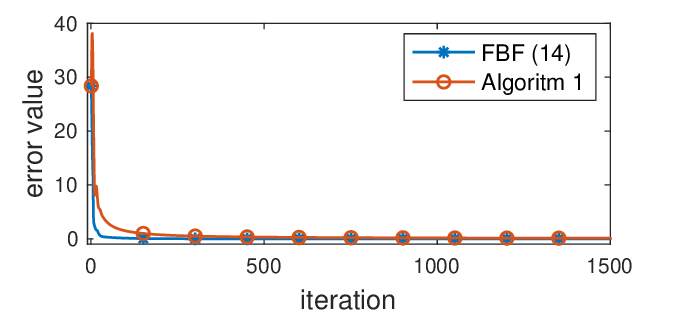}
  \vspace*{-0.1cm}
  \renewcommand{\thefigure}{2}
  \caption{Approximation error $\|(P_{\overline{B}(0, r)} \circ \Talpha)(\bm{\xi}_n) - \bm{\xi}_n\|_2$ of variational equilibria for \eqref{eq:Cournot-Nash}}
  \label{fig:resi}
\end{figure}
\subsection{Linearly-coupled aggregative game}
We apply Algorithm 1 to the HNEP in a scenario of the linearly-coupled aggregative game \cite[Section IV.B]{g.belgioiosoSemiDecentralizedGeneralizedNash2023a} which is found in many engineering problems, e.g., the electric vehicle charging control problem \cite{deoriPriceAnarchyElectric2018} and the energy management problem for peer-to-peer trading systems \cite{g.belgioiosoEnergyManagementPeertopeer2020}. 

As a special instance of Problem \ref{prob:GNEP}, the linearly-coupled aggregative game is given by 
\begin{align}
  \label{eq:Cournot-Nash}
  \begin{split}
    &\find \bar{\x} = (\bar{x}_1, \dots, \bar{x}_m) \in \textstyle \bigtimes_{i\in \mathcal{I}} \R^{M} = \R^{mM} \ \st (\forall i \in \mathcal{I})\\
    &\bar{x}_i \in \argmin[x_i \in \R^M] \iota_{C_i}(x_i) + \f_i(x_i; \bar{\x}_{\smallsetminus i}) + (\iota_{D}\circ\Li)(x_i; \bar{\x}_{\smallsetminus i}),
  \end{split}
\end{align}
where $M\in \N\setminus\{0\}$ and for every $i \in \mathcal{I}$, 
the player $i$'s local feasible set and $i$'s cost function are given respectively by 
\begin{align}
  C_i &:= \bigtimes_{j \in \{1, \dots, M\}}[a_{i, j}, b_{i, j}] \subset \R^M, \\
  \f_i(\x) &:= \left(\frac{1}{m} \sum_{i\in I} \bm{W}x_i - p\right)^Tx_i, 
\end{align}
with $a_{i, j} < b_{i, j}$ for every $j \in \{1, \dots, M\}$, 
$p \in \R^M_{++}$, and a nonnegative diagonal matrix $\bm{W}\in\R^{M\times M}$. 
In the problem \eqref{eq:Cournot-Nash}, $\Li$ and $D$ are given respectively by $\Li:\R^{mM}\rightarrow\R^M:\x:=(x_1, \dots, x_m) \mapsto \sum_{i\in\mathcal{I}} x_i$ and $D := \{y\in\R^M\mid y \leq c\}$\footnote{ $y \leq c$ means that $y_j \leq c_j$ for each component $j\in\{1, \dots, M\}$.}, where each component of $c\in\R^M_{++}$ represents the upper bound of each component of $\Li(\x)\in\R^M$. 
We can check that the problem \eqref{eq:Cournot-Nash} enjoys 
the setting of Problem \ref{prob:GNEP} and Proposition \ref{prop:zero_closed_convex}. 

A variational equilibrium, but an anonymous one, for the problem \eqref{eq:Cournot-Nash} can be approximated iteratively \cite{tsengModifiedForwardBackwardSplitting2000}, \cite[Proposition 17.10]{yamadaMinimizingMoreauEnvelope2011} by a slight modification\footnote{Although a variational equilibrium can be estimated by the simple algorithm $\bm{\xi}_{n+1}=\Talpha(\bm{\xi}_n)$, we used the algorithm in \eqref{eq:ex_FBF} with the projection $P_{\overline{B}(0, r)}$ for a fair comparison to the proposed algorithm.} of the forward-backward-forward splitting (FBF) algorithm:
\begin{align}
  \label{eq:ex_FBF}
  (n \in \N) \ \bm{\xi}_{n+1} = (P_{\overline{B}(0, r)} \circ \Talpha)(\bm{\xi}_n),
\end{align}
where $\Talpha:\R^{mM} \times \R^{M} \rightarrow \R^{mM} \times \R^{M}$ 
is defined as \eqref{eq:Talpha} and 
$\overline{B}(0, r) \subset \R^{mM} \times \R^{M}$ is a sufficiently large closed ball (see also Remark \ref{rem:ass} (ii)). 

\subsection{Designing each player's upper-level cost function}
Inspired by \cite[Sec. VI.A]{10018595}, for every $i \in \mathcal{I}$, we define the player $i$'s upper-level cost function $\uf{i}$ in the HNEP \eqref{eq:HierarchicalGame_0} by 
\begin{equation}
  \label{eq:evaluation}
  \uf{i}(\x) := \frac{1}{2} \biggl( 
    \|x_i - t_i\|^2 + 
    \sum_{j \in \mathcal{I} \setminus \{i\}} \|x_i - x_j\|^2 \biggr), 
\end{equation}
where $t_i \in \R^M$ is the player $i$'s target point, 
the second term is the sum of the square distances between the player $i$'s strategy and the other players' strategies. 
Note that the player $i$'s local target point $t_i$ can be determined according to the player $i$'s hope. 
Then, we consider the following problem for finding a special solution of the HNEP \eqref{eq:HierarchicalGame_0}:
\begin{align}
  \begin{split}
    \label{eq:VI_HierarchicalGame_F_B_app}
    &\find \bm{\xi}^{\star} \in \Fix(\Talpha) \ \st \\
    &(\forall \bm{\zeta} \in \Fix(\Talpha)) \ \product[]{\uG_{\FBF}(\bm{\xi}^{\star})}{\bm{\zeta} - \bm{\xi}^{\star}}\geq 0,
  \end{split}
\end{align}
where $\uG_{\FBF}:\R^{mM} \times \R^{M}\rightarrow\R^{mM} \times \R^{M}:\bm{\xi}:=(\x, u) \mapsto (\uG(\x), \bm{0}_{\R^{mM}})$ and 
$\uG : \R^{mM}\rightarrow\R^{mM}$ is defined as \eqref{eq:uG}. 
In this case, Assumption \ref{asm:paramonotone} (i) is satisfied. 
As mentioned in Remark \ref{rem:ass} (ii), 
we can find a solution of the problem \eqref{eq:VI_HierarchicalGame_F_B_app} 
by a slight modification of \eqref{eq:HSDM} (Algorithm \ref{alg:HSDM}):
\begin{equation}
  \label{eq:HSDM_app}
  \begin{split}
    (n\in\N) \ \bm{\xi}_{n+1} = &(P_{\overline{B}(0, r)} \circ \Talpha)(\bm{\xi}_n) - \\
    &\lambda_{n+1} \uG_{\FBF}\bigl((P_{\overline{B}(0, r)} \circ \Talpha)(\bm{\xi}_n)\bigr). 
  \end{split}
\end{equation}

\subsection{Numerical experiments}
We compared a solution of the problem \eqref{eq:VI_HierarchicalGame_F_B_app} with a variational equilibrium for the problem \eqref{eq:Cournot-Nash}. 
For each $i \in \mathcal{I}$ and $j \in \{1, \dots, M\}$, the parameters\footnote{We set the parameters in the problem \eqref{eq:Cournot-Nash} along the setting found in \cite{e.benenatiOptimalSelectionGeneralized06} where an equilibrium selection over $\bm{\mathcal{V}}$ is considered as a convex optimization based on a hybrid steepest descent method \cite{oguraNonstrictlyConvexMinimization2003}.} were set as follows: 
$(m, M) = (6, 3)$, 
$b_{i, j} = 100$, 
$c = (120, 120, 120)^T$, 
$a_{i, j} \in [-1, 1]$, 
$p \in [0, 10]\times[0, 10]\times[0, 10]$, 
each nonnegative entry of $\bm{W}$ belongs to $[0, 1]$, 
and $t_i \in C_i$. 
In this experiment, we applied the FBF algorithm \eqref{eq:ex_FBF} 
to the problem \eqref{eq:Cournot-Nash} and the proposed algorithm \eqref{eq:HSDM_app} to the problem \eqref{eq:VI_HierarchicalGame_F_B_app}, 
where the parameters in \eqref{eq:FBF}, \eqref{eq:ex_FBF}, and \eqref{eq:HSDM_app} were set as $\gamma = 0.25$, $\alpha = 0.75$, $r = 10^{15}$, and $\lambda_n = 1/(n+3)$\footnote{
  We employed the stepsize $\lambda_n:=(1/(n+3))$ for the proposed algorithm \eqref{eq:HSDM_app} to enjoy the numerical stability. 
  We note that such a $\lambda_n$ satisfies the conditions (H1) $\lim_{n\to\infty}\lambda_n = 0$ and (H2) $\sum_{n\in\N}\lambda_n=\infty$, implying that the proposed algorithm \eqref{eq:HSDM_app} has a guarantee of convergence in the sense of Theorem \ref{thm:HSDM_convergence}.}. 
The initial point $\bm{\xi}_0 = (\x_0, u_0)$ was randomly generated from uniform distribution. 

Fig. \ref{fig:resi} shows an approximation error 
$\|(P_{\overline{B}(0, r)} \circ \Talpha)(\bm{\xi}_n) - \bm{\xi}_n\|_2$ 
between the $n$th estimate $\bm{\xi}_n$ and variational equilibria for the problem \eqref{eq:Cournot-Nash}. 
From Fig. \ref{fig:resi}, we can see that both algorithms approximate a variational equilibrium for the problem \eqref{eq:Cournot-Nash} 
(Note: each element of $\Fix(P_{\overline{B}(0, r)} \circ \Talpha)$ is a variational equilibrium for the problem \eqref{eq:Cournot-Nash}). 
Fig. \ref{fig:cost_u} shows the value of each player's upper-level cost function $\uf{i} \ (i\in\mathcal{I})$ defined in \eqref{eq:evaluation}. 
From Fig. \ref{fig:cost_u}, we can see that for every player $i\in \mathcal{I}$, the values of $\uf{i}$ by Algorithm \ref{alg:HSDM} are lower than the values of $\uf{i}$ by \eqref{eq:ex_FBF}.
To sum up, as we expected, the proposed algorithm \eqref{eq:HSDM_app} can iteratively approximate a variational equilibrium, of the problem \eqref{eq:VI_HierarchicalGame_F_B_app}, which also decreases the values of each player's upper-level cost function $\uf{i}$. 
\section{Conclusion}
We proposed a new equilibrium selection method, 
without assuming any trusted center or any randomness assumption, 
achievable by solving the hierarchical Nash equilibrium problem (HNEP). 
We also proposed an iterative algorithm for the HNEP as an application of the hybrid 
steepest descent method to the variational inequality defined over the 
fixed point set of a quasi-nonexpansive operator. 
Numerical experiments showed the effectiveness of the proposed equilibrium selection method. 

\bibliographystyle{IEEEbib}
\bibliography{./bib/refs, ./bib/strings, ./bib/IEEEabrv}

\end{document}